\documentstyle[12pt]{article}
\def\halm{\vrule height1ex width.9ex depth.1ex} 
\def\Ker{{\rm Ker\;}}
\def\Coker{{\rm Coker\;}}
\def\Im{{\rm Im\;}}
\def\H{{\rm H}}
\def\db #1 {{\bf#1}}

\def\E{{\rm E}}
\def\A{{\rm A}}
\def\III {{\rm III}}
\begin{document}
\title{Weak approximation and Manin group of R-equivalences}
\author{Nguy\^e\~n Qu\^o\'{c} Th\v{a}\'ng\thanks
{Support by Department of Mathematics and Statistics, McMaster University is
greatly acknowledged}}
\date{}
\maketitle
\begin{center}Department Mathematics and Statistics
\\
McMaster
University,
Hamilton
\\
Ontario, Canada L8S 4K1
\end{center}
\begin{abstract} We extend one exact sequence of Colliot-Th\'el\`ene and Sansuc
for tori over number fields to one for arbitrary connected groups.\footnote{
Mathematics Subject Classification : Primary 20G 10, Secondary 11E72, 18G50}
\end{abstract}

{\bf Introduction.}
Let $G$ be a connected linear algebraic  group defined over a number field $k$.
Denote by $$\A(G) = \prod _v G(k_v)/Cl (G(k))$$ the obstruction
to weak approximation in $G$, $$\III (G) = \Ker (\H^1(k,G) \to \prod_v
\H^1(k_v,G))$$ the Tate - Shafarevich group of $G$ where $\H^1(.,.)$
denotes Galois cohomology and $k_v$ denotes the completion of the field $k$
at a valuation $v$. We denote also
by $G(k)/R$ the  Manin  group of $R$-equivalences
of $G(k)$, $RG(k)$ the subgroup of $G(k)$ consisting of all elements
of $G(k)$ which are $R$-equivalent to 1 in $G(k)$ (see [M]). One defines
similar
groups in the case of local fields $k_v$ which are completions  of $k$.
\\
If $G=T$ is a $k$-torus, then Voskresenski [V1] showed that there is an exact
sequence connecting arithmetic, geometric and cohomological invariants
of $T$
$$ 1 \to \A(T) \to \H^1(k,Pic\bar {X})^* \to \III (T) \to 1, \eqno{(1)}$$
where $X$ is a smooth compactification of $T$ over $k$,
$\bar X = X \times_k \bar k$ and $()^*$ means taking the character group
with values in ${\bf Q}/{\bf Z}$.  Later on, Sansuc [Sa] showed that this
sequence also holds for arbitrary connected $k$-group $G$ over  a number field
$k$ without
factors of type $\E_8$, thus also for any connected $k$-group by
combining this with the result on the Hasse principle for $\E_8$ by Chernousov.
In their fundamental paper [CTS] Colliot-Th\'el\`ene and Sansuc  have
established, among the others,
the following exact sequence which connects various
important arithmetic and geometric invariants of algebraic tori defined
over a number field $k$. Namely we have
$$ 1 \to \III (S) \to T(k)/R \stackrel{\beta}{\to}
\prod_v T(k_v)/R \to \A(T) \to 1, \eqno{(2)}$$
an exact sequence where  $S$ is the Neron - Severi torus of $T$ which comes
froma flasque resolution of $T$.
\\
This exact sequence gives us the information about
the kernel of $\beta$ : it is the group $\III (S)$, and the cokernel
of $\beta$ : it is the group $\A(T)$. It is a natural question to ask if one
can extend (2) to the case of connected reductive groups over global fields.
Denote by $\III RG$ the "local-global"
group of $R$-equivalent classes modulo $RG(k)$
which are trivial locally everywhere, i.e.,
$$\III RG = \cap_v RG(k_v)/RG(k).$$
Then (2) can be written in the following form
$$ 1 \to \III RT \to T(k)/R \to \prod_v T(k_v)/R \to \A(T) \to 1. \eqno{(3)}$$
We show in this note that (3) also holds for any connected reductive $k$-group
$G$.
Namely we have
\\
\\{\bf Theorem.} {\it  For any connected group G defined over a global field k,
assumed reductive if $char.k >0$, there is an exact sequence of groups}
$$ 1 \to \III RG \to G(k)/R \to \prod_v G(k_v)/R \to \A(G) \to 1.$$
\\
\\
We give some preliminary results in Section 1 and we prove the main result
in Section 2. In Section 3 we give a cohomological interpretation of the
group $\III RG$. I would like to thank J.-L.Colliot-Th\'el\`ene for the
interest
in the paper.
\section {Some auxialary results}
We need the following auxialary results in the arithmetic of algebraic groups.
\\
\\
{\bf 1.1. Theorem.}
{\it Any semisimple simply connected group defined
over a global field k satisfies weak approximation.}
\\
\\
This result is due to Kneser and Harder.
We refer to [Sa] for  historical remarks and references.
\\
The following lemma of Langlands allows one to reduce many problems
in reductive groups to the case of a reductive group with simply connected
semisimple part.
\\
\\
{\bf 1.2. Lemma.} [L] {\it  Let G be a connected reductive group
defined over a field k. Then there exists a connected reductive group
H, an induced torus Z and a central extension}
$$ 1 \to Z \to H \to G \to 1,$$
{\it all defined over k.}
\\
\\
In literature, such central extension of $G$ is called a $z$-extension of
$G$.
\\
\\
{\bf 1.3. Theorem.} {\it Let
G be a connected reductive group defined over a local or global field k.
Then the Manin group $G(k)/R$ is finite.}
\\
\\
In [G2] Gille proved this theorem in the case of number field, but
the proof can be modified to give other cases as well (see [T1,T2]).
\\
\\
{\bf 1.4. Proposition.} {\it Let G be a unirational group over a field
k and S a finite set of discrete valuations of k. Then via the diagonal
embedding, the closure $Cl(G(k))$ of $G(k)$
in $\prod_{v \in S}G(k_v)$ in the product topology contains an open
subgroup of the latter.}
\\
\\
This result is implicitly
 due to Kneser, who used it in the study of
strong approximation. The proof follows from the Implicit Function Theorem.
\\
\\
{\bf 1.5. Theorem.} [CTS] {\it Let T be a torus defined over a field k. Then
there exists a flasque resolution of T, i.e., an exact sequence of $k$-tori}
$$ 1 \to S \to N \to T \to 1,$$
{\it where N is an induced torus and S is a flasque torus. Moreover
$T(k)/R \simeq \H^1(k,S)$.}
\section {Main theorem and the proof}
First we need the following
\\
\\
{\bf 2.1. Lemma.} {\it Let G be a connected reductive and H a z-extension of
G, all defined over a field k. Assume that either $H(k)/R$ or $G(k)/R$
is finite. Then there are canonical isomorphism
of groups and bijection of factor sets, respectively}
$$G(k)/R \simeq H(k)/R,~ \A(G,S) \simeq \A(H,S),$$
{\it where $\A(G,S)$ denotes the factor $\prod_{v \in S} G(k_v)/Cl(G(k))$
with S a finite set of discrete valuations of k and $Cl(G(k))$ the closure
of $G(k)$ via diagonal embedding  in the product topology of $\prod_{v \in S}
G(k_v)$.}
\\
\\
{\it Proof.} Let $1 \to Z \to H\stackrel{\pi}{\to} G \to 1$
be the corresponding $z$-extension.
Since $Z$ is an induced torus, it has trivial cohomology so $\pi$ is
surjective while restricted to $H(K)$ for any field extension $K$ of $k$.
It is clear that $\pi(RH(k)) \subset RG(k)$ hence $\pi$ defines a surjective
homomorphism of groups $$\pi' : H(k)/R \to G(k)/R.$$
It is well-known that $H \simeq Z \times G$ (birationally equivalent) and
we know  by [CTS] that this induces a bijection
$$H(k)/R \leftrightarrow G(k)/R.$$ Since $H(k)/R$ is finite, $\pi'$ is an
isomorphism of groups.
\\
To prove the second bijection, it suffices to prove the following equality
$$ \pi ^{-1} (Cl(G(k)) = Cl(H(k)).$$
One inclusion is obvious. Let $ h \in \pi ^{-1}(Cl(G(k)).$ Then
$\pi(h) = lim_n g_n,$ with $g_n \in G(k)$. There are $h_n \in H(k)$
such that $g_n = \pi(h_n)$, so $lim_n \pi(h_n^{-1}h) = 1$.
Let $U_n$ be a system of open neibourghoods of 1 in $G_S =
\prod_{v \in S}G(k_v)$ such that $U_n \subset U_{n-1}$, $\cap_n U_n =\{1\}$
and $\pi(h_n^{-1}h) \in U_n$. Thus for all $n$ we have
$$ \pi^{-1}(1) = \pi^{-1}(\cap_n U_n) = \cap_n \pi^{-1}(U_n) = \prod_{v \in
S}Z(k_v),$$
and
$$ h_n^{-1}h \in \pi^{-1}(U_n) \eqno{(4)}$$
It is well-known by a result of Grothendieck - Rosenlicht
that connected reductive groups are unirational over the field of definition.
Thus by Proposition 1.4,
$Cl(H(k))$ is an open subgroup of $H_S$, hence $Cl(H(k)) \prod_{v \in S}Z(k_v)$
is an open neighbourhood of $\prod_{v \in S} Z(k_v)$ in $\prod_{v \in
S}H(k_v)$.
Therefore for $n$
large we have
$$ \pi^{-1}(U_n) \subset Cl(H(k)) \prod_{v \in S}Z(k_v) \subset Cl(H(k)),$$
since $Z$ has weak approximation property.
Thus for $n$ large we derive from  (4)
$$ h \in h_n \pi^{-1}(U_n) \subset h_n Cl(H(k))=Cl(H(k))$$
and the lemma follows. \halm
\\
\\
{\bf 2.2. Lemma.} {\it Let G be a connected reductive
group defined over a field  $k$, H a z-extension of G. For a finite set
S of valuations of k, let $RG_S = \prod _{v \in S} RG(k_v)$. }
\\
$a)$ {\it If $RG_S \subset Cl(G(k))$ then $RH_S \subset Cl(H(k))$.}
\\
$b)$ {\it Conversely, assume that $H(k_v)/RH(k_v)$ is finite for $v \in S$
and $RH_S \subset Cl(H(k))$. Then $RG_S \subset Cl(G(k))$.
}
\\
\\
{\it Proof.} $a)$ We have $\pi (RH_s) \subset RG_S$ where we denote by
the
same symbol $\pi$ the homomorphism $H_S \to G_S$ induced from $\pi : H \to G$.
Hence
$$RH_S \subset \pi^{-1}(RG_S) \subset
\pi^{-1}(Cl(G(k))) = Cl(H(k))$$ by Lemma 2.1.\\
$b)$ By Lemma 2.1 we have $\pi(RH_S) = RG_S$, thus
$$RG_S \subset \pi(Cl(H(k))) = Cl(G(k)).$$
The lemma is proved. \halm
\\
\\
The following was first mentioned in [V2] in the case the finite set of
valuations consists of only one element.
\\
\\
{\bf 2.3. Lemma.} {\it Let T be a torus defined over a field k, S a finite
set of valuations of k. Then $RT_S \subset Cl(T(k))$.}
\\
\\
{\it Proof.} Consider the following flasque resolution of $T$
$$ 1 \to S \to N \stackrel{p}{\to} T \to 1.$$
Since this is also a flasque resolution of $T$ over $k_v$, from the proof
of Theorem 1.5 (in [CTS]) it follows that
$$ p(N(k)) = RT(k),~ p(N(k_v)) = RT(k_v)$$
for any $v \in S$ and since $N$ has weak approximation over $k$, it follows
that $$ \prod_{v \in S} RT(k_v)  = p(\prod_{v \in S}N(k_v))
= p(Cl(N(k)))  \subset Cl(T(k)).$$
The lemma is proved. \halm
\\
\\
For arbitrary $k$-group it is not known if the Lemma 2.3 is true in
general,
but in the important case of connected reductive groups over
global fields we have the following result, which answers a (implicit) question
of
Voskresenski [V2] in this case.
\\
\\
{\bf 2.4. Theorem.} {\it Let k be a global field and G a connected reductive
group over k. For any finite set S of valuations of k we have }
$$ RG_S \subset Cl(G(k)).$$
{\it Proof.} By 2.3 we may assume that $G$ is not a torus.
By 1.3 and 2.2 we may assume that the semisimple part $G'$ of $G$ is simply
connected.  By enlarging $S$ we may also assume that $S$ contains all
archimedean valuations of $k$. Let $G = G' T$, where $T$ is a central torus.
We have the following exact sequence of $k$-groups
$$ 1 \to G' \to G \stackrel{p}{\to} T' \to 1,$$
where $T'$ is a factor of $T$. Consider the following commutative diagram,
where
the rows are exact \\
\\
\noindent
\hspace*{3.5cm}$G(k) \to~ T'(k)~ \stackrel{\delta}{\to}~ \H^1(k,G')$\\
\\
\hspace*{3.9cm}$\downarrow \hspace{1.4cm}\downarrow
\hspace{1.5cm}\downarrow \gamma$\\
\\
\hspace*{3.5cm}$G(S) \to T'(S) \stackrel{\delta_S}{\to} \prod_{v \in S}
\H^1(k_v,G'),$\\
\\
where $G(S) = \prod_{v \in S}G(k_v)$, $T'(S) = \prod_{v \in S} T'(k_v)$.
Since $G'$ is simply connected, we know that $\gamma$ is bijective.
Let $g \in RG_S$. Then $p(g) \in RT_S \subset Cl(T'(k))$ by 2.3. Therefore
there is a sequence $(t_n)$, $t_n \in T'(k)$ such that
$$ p(g) = lim_n t_n,$$
hence $$\delta (p(g)) = \delta(lim_n t_n) = 1.$$
Since the product $\prod_{v \in S} \H^1(k_v,G')$ is finite, for large $n$ we
have
$$ \delta (t_n) = 1,$$
i.e.,  for large $n$, there are $g_n \in G(k)$ such that $t_n = p(g_n)$.
Therefore
$$lim_np(g_n^{-1}g) = 1.$$
By the same argument we used in the proof of 2.1, and by using weak
approximation
property in simply connected groups (see 1.1) we see that for large $n$
$$ g_n^{-1} g \in G'(S) Cl(G(k)) \subset Cl(G(k)),$$
i.e., $g \in Cl(G(k))$. \halm
\\
\\
To deduce the relation between weak approximation and R-equivalence we need the
following well-known
\\
\\
{\bf 2.5. Lemma.} {\it With the same notation as in 2.1, $RG_S$ is an open
subgroup of $\prod_{v \in S} G(k_v)$.}
\\
\\
{\it Proof.} It follows from the fact that $G$ is unirational and that affine
spaces
have only trivial $R$-equivalent classes. \halm
\\
\\
Now  we have
\\
\\
{\bf 2.6. Theorem.}
{\it Let $k$ be a global field and G a connected k-group, which is reductive if
$ char.k >0$. For any finite set S of valuations
of k  we have a bijection}
$$ \A(S,G) \leftrightarrow \Coker (G(k)/R \to \prod_{v \in S} G(k_v)/R).$$
{\it In particular, the obstruction $\A(S,G)$ has a natural group structure
inherited
from that of $\Coker (G(k)/R \to \prod_{v \in S} G(k_v)/R)$.}
\\
\\
{\it Proof.}
By 2.4 we have $$RG_S \subset Cl(G(k)),$$
and by 2.5 $$Cl(G(k)) = RG_S G(k).$$ Thus
$$\Coker (G(k)/R \to \prod_{v \in S} G(k_v)/R ) =
\prod_{v \in S}G(k_v)/Cl(G(k)) = \A(S,G).$$
The theorem is proved. \halm
\\
\\
It is well-known (see [CTS], [Sa]) that for almost all $v$,
$G(k_v)/R = 1$ and there is a finite set $S_0$ of valuations
of $k$ such that $G$ has weak approximation over $k$ with respect to any finite
set of valuations of $k$ outside $S_0$. Thus we have \\
\\
{\bf 2.7. Theorem.} {\it Let the notation be as in 2.6 and the introduction.
Then we have
the following
exact sequence of groups}
$$ 1 \to \III RG \to G(k)/R \to \prod_{v} G(k_v)/R \to \A(G) \to 1.
\eqno{(5)}$$
\section {A cohomological interpretation of $\III RG$.}
As we have seen, in the case  $G=T$ is a
torus, the group $\III RG$ is the Tate - Shafarevich
group of the Neron - Severi torus $S$ of $T$.
We would like to find a similar interpretation for the case of reductive groups
over number fields.
\\
{}From [Sa] we know that for any connected reductive group $G$ over a field
$k$,
there is a number $n$, induced $k$-tori $P$, $Q$ and a central $k$-isogeny
$$ 1 \to F \to \tilde G \times P \to G^n \times Q \to 1,$$
where $\tilde G$ is a semisimple simply connected $k$-group (namely the simply
connected
covering of  $G'^n$, $G'$ being the semisimple part of $G$.
{}From the proof in the case of tori and the case of reductive groups, we see
that for
$n$ smallest possible, the groups $P$, $Q$ are determined uniquely. The
corresponding
group $\tilde G \times P$ is called {\it canonical special covering of G} and
$F$ the {\it canonical special fundamental group of G}. In the case $G$ is
semisimple,
these notions are just the usual notions of
simply connected covering  and fundamental group.
\\
First we consider the case $G=T$ is a torus. There is a canonical special
covering
of $T$
$$ 1 \to F \to P \to T^n \times Q \to 1,$$
where $P, Q$ are induced tori, $F$ a finite group.
\\
For any extension $K$ of $k$ we have the following exact sequence  of groups
$$ P(K) \stackrel{p}{\to} (T^n \times Q)(K) \stackrel{\delta}{\to} \H^1(K,F)
\to 0,$$
and similar sequence when we take the group of R-equivalences
$$ 1 \to T^n(K)/R \to \H^1(K,F)/R \to 0,$$
where in the case of $\H^1(K,F)$,
we consider the $R$-equivalence induced from $(T^n \times Q)(K)$, (see [G1]).
Therefore
$$T^n(K)/R = \H^1(K,F)/R,$$
and we have
$$\III RT^n = \III F/R := \Ker (\H^1(k,F)/R \to \prod_v \H^1(k_v,F)/R).$$
Thus from (5) we deduce
the following exact sequence for tori
$$ 1 \to \III F/R \to T^n(k)/R \to \prod_v T^n(k_v)/R \to \A(T^n) \to 1.
\eqno{(6)}$$
Next we discuss the analog of (6) in the case of connected reductive groups.
\\
Let $G$ be a connected reductive $k$-group with canonical special fundamental
group
$F$. There is an embedding
$$ \H^1(k,F) \to \H^1(k(t),F),$$
where $k(t)$ is the rational function field in the variable $t$. There is an
equivalence
relation
$R$ on $\H^1(k,F)$, defined as follows :
\\
for $x, y \in \H^1(k,F)$, $x \sim_R y$
iff there exists $z(t) \in \H^1(k(t),F)$ such that there are specializations
$t \mapsto 0$, $t \mapsto 1$ with $z(0) =x, z(1)=y.$
\\
We denote by $\H^1(k,F)/R$
the set of $R$-equivalence classes.
Let $\delta : (G^n \times Q)(k) \to \H^1(k,F)$ be the connecting map
and we denote by $\alpha$ the composite map
$$ G^n(k)/R \to \Im(\delta)/R \to \H^1(k,F)/R.$$
The second map is induced from the embedding $\Im(\delta) \to \H^1(k,F)$.
We have
\\
\\
{\bf 3.1. Theorem.} {\it Let G be a connected reductive group defined over a
global
field k and other notation be as above.}
\\
$a)$ {\it If k is a global field of characteristic $>0$ and G contains no
anisotropic factors of type $^2\A_n$, then the following sequence is exact}
$$ 1 \to \III F/R \to G^n(k)/R \to \prod_v G^n(k_v)/R \to \A(G^n) \to 1.
\eqno{(7)}$$
$b)$ {\it In general, we have the following exact sequence of groups}
$$ 1 \to \alpha^{-1}(\III F/R) \to G^n(k)/R \to \prod_v G^n(k_v)/R \to \A(G^n)
\to 1.$$
\\
{\it Proof.} $a)$ By a well-known theorem of Kneser [K] (in characteristic 0)
and
Bruhat - Tits [BT] (in characteristic $>0$), $\H^1(k_v,\tilde G) =0$ where
$\tilde G$ is a semisimple
simply connected group  and $k_v$ is a non-archimedean local field. Also, if
$k$ is a global field of characteristic $>0$, $\H^1(k,\tilde G)=0$ for simply
connected
semisimple group $\tilde G$ by [Ha]. It is well-known that $\tilde G(k)$ is
projectively
simple, hence also $\tilde G(k)/R=1$
for simply connected $k$-group $\tilde G$ containing no anisotropic factors of
type
$^2\A_n$ (see [PR] and [Ti2]).
Thus the  exact sequence (7) follows as in the case of tori.\\
$b)$ In the general case, for any non-archimedean valuation $v$ we
still have $G^n(k_v)/R \simeq \H^1(k_v,F)/R$. It is known and
easy to prove that $H(\db R )/R=1$ for any connected
\db R -group. Hence
\begin{tabbing}
$\Ker (G^n(k)/R \to \prod_v G^n(k_v)/R)$ \= $=\Ker (G^n(k)/R \to \prod_v
\H^1(k_v,F)/R)$
\\
\\
\>$= \alpha^{-1}( \III F/R).$
\end{tabbing}
The theorem is proved. \halm
\\
\\
Conjecturally, $\alpha$ is an isomorphism and $\tilde G(k)/R=1$ for any simply
connected
semisimple group  $\tilde G$ over a
global field $k$. Moreover, the term $\III F/R$ might be
expressed  as a function of $S:= Pic (V(G))^*$, where $k$ is assumed a number
field
and $V(G)$ is a smooth compactification of $G$.

\end{document}